\documentclass[12pt]{article}
\usepackage{tikz}
\usepackage{amsmath}
\usepackage{amssymb}
\usepackage{graphicx}
\usepackage[numbers,sort&compress]{natbib}
\usepackage[colorlinks,
            anchorcolor=blue,
            citecolor=green]
            {hyperref}
\title{Maps from the Dbar-data to the derivative NLS potentials}

\author{Junyi Zhu$^a$,
Xueru wang$^b$ \thanks{E-mail: xr\_wangmath@163.com}, Huan Liu$^a$
\\
\scriptsize{$^a$ School of Mathematics and Statistics, Zhengzhou University, Zhengzhou, Henan 450001, China}\\
\scriptsize{$^b$ School of Mathematics and Statistics, Xinyang Normal University, Xinyang, Henan 464000, China}
}

\date{}
\begin{document}
\newtheorem{theorem}{Theorem}
\newtheorem{definition}{Definition}
\newtheorem{lemma}{Lemma}
\newtheorem{corollary}{Corollary}
\newtheorem{proposition}{Proposition}
\newtheorem{remark}{Remark}
\hypersetup{CJKbookmarks=true}

\maketitle
\begin{abstract}
We develop a decomposition technique to investigate maps from the Dbar data to the derivative NLS potentials. We define a new integral operator which is equivalent to a modified Dbar problem relevant to the derivative NLS spectral problem, including the Kaup-Newell spectral problem, the Chen-Lee-Liu spectral problem and the Gerdjikov-Ivanov spectral problem. The well-posedness for the integral operator with the kernel $\mathrm{e}^{ik^2x\sigma_3}$ is discussed via the decomposition technique. Under different constraint conditions, the derivative NLS potentials are reconstructed via the Dbar data. We discuss the symmetry condition of the Dbar data as well as the compatibility of the derivative NLS spectral problem with the modified Dbar problem. We show that the maps from the Dbar data to the derivative NLS potentials are Lipschitz continuous.

\end{abstract}

\section{Introduction}
The derivative NLS (DNLS) equation
\begin{equation}\label{dnls3}
iu_t+u_{xx}+i\varepsilon(|u|^2u)_x=0, \quad \varepsilon=\pm1,
\end{equation}
is a canonical dispersive equation derived from parallel propagation of nonlinear low-frequency waves in high-$\beta$
plasma \cite{pf14-2733}. This equation is integrable \cite{jmp19-798} and has many important applications in different fields, such as long-wavelength dynamics of dispersive Alfv\'en waves in the plasma physics field \cite{pf14-2733,jpp16-321,jpsj41-265}, the subpicosecond or femtosecond pulse in a single-mode fiber \cite{ol12-516}, the propagation of nonlinear pulses in optical fibers \cite{AG2001}, and so on.

The DNLS equations have several different forms. The other well-known forms are the Chen-Lee-Liu equation \cite{ps20-490}
\begin{equation}\label{CLL1}
iu_t+u_{xx}+i\varepsilon |u|^2u_x=0, \quad \varepsilon=\pm1,
\end{equation}
and the Gerdjikov-Ivanov equation \cite{bjp10-130}
\begin{equation}\label{GI1}
iu_t+u_{xx}-i\epsilon u^2\bar{u}_x+\frac{1}{2}\epsilon u^3{\bar{u}}^2=0,\quad \epsilon=\pm1,
\end{equation}
where $\bar{u}$ denotes the complex conjugate of $u$.
We note that under some gauge transformations, DNLS equation \eqref{dnls3} can be mapped into the Chen-Lee-Liu equation \eqref{CLL1} and the the Gerdjikov-Ivanov equation \eqref{GI1}.

The $\bar\partial$ (Dbar) problem can be regarded as a generalization of the Riemann-Hilbert problem. Since Beals and Coifman introduced the Dbar problem to study the integrable system \cite{B-C1981,B-C1982,ip5-87}, Dbar problem has played an important role to construct explicit solutions \cite{ip4-123,D-L2007,faa19-89,jpa21-L537,KBG1992,KBG1993,Santini2003,sam69-135,jpa46-035204,mpag17-49,aml66-47,wm115-103051,sam148-433,amp13-51,pla499-129359,jgp177-104550,aml107-107297,
aml117-107143,FEGE2022,jmp62-093510,jnmp28-492,aml134-108378,nd111-3689,nd111-5655,nd111-8673,jnmp30-201,aml140-108589,prsa481-20250342,ps100-115237,prsa481-20240764}, and to discuss the asymptotic behaviors of the nonlinear integrable system \cite{imrn2006-48673,jmp63-113504,cmp402-2879,jde329-31,aim426-109088,jde386-214,jde447-113685,qtds24-106,pd456-133915,jde377-121,jmp66-041502,aim484-110552,sam155-e70093}.
A prior estimates on integral operator to the nonhomogeneous Dbar equation without physical variables ($x,y,t$) have been well discussed \cite{INV1962}.  A matrix integrable operator $\mathcal{K}$ on a bounded domain $\mathcal{D}$ is considered in \cite{non37-085008}, and the resolvent of the integral operator $Id - \mathcal{K}$ is discussed. A decomposition technique has been developed to discuss the Dbar well-posedness relevant to the AKNS spectral problem \cite{2603.18930v2}.

If the Dbar problem is applied to the DNLS spectral problems, the matrix exponentials $\mathrm{e}^{\pm ik^2 x}$ depending on physical variable $x$ will be introduced into the spectral transformation matrix of the Dbar problem. And the Dbar problem is equivalent to an integral equation, with exponentials $\mathrm{e}^{\pm ik^2 x}$ in the kernel. We note that the integral may not convergent for $\mp x{\rm Im}(k^2)>0$ as $k\in \mathbb{C}$ and $x\in \mathbb{R}$.
Here, we will extend the decomposition technique to investigate the spectral problems of the DNLS equations.
To this end, we introduce a new integral operator defined by $\psi RT_\mathbb{C}$ and show that the new integral equation $\psi=\psi_0+\psi RT_\mathbb{C}$ is equivalent to a modified Dbar problem with normalization condition $\psi(x,k)\to\psi_0(x), k\to\infty$. We show that the inverse operator $(\mathcal{I}-RT_\mathbb{C})^{-1}$ exists and is bounded, then the function $\psi$ was exist uniquely and given by
\begin{equation}\label{kna7}
\psi(k;x)=I(\mathcal{I}-RT_\mathbb{C})^{-1}(k;x).
\end{equation}
We extend the Dbar dressing method, and prove that the modified Dbar problem is related to the DNLS spectral problem.
Compared with the technique used in \cite{2603.18930v2}, in the prior estimates, the convergent domain of the associated integral operator will be more complicated, while the present extended decomposition technique is simpler and more straightforward.

The DNLS equations, including the Chen-Lee-Liu equation and the Gerdjikov-Ivanov equation, are important integrable equations,  and can be solved by the inverse scattering transform method based on different linear spectral problems \cite{jmp19-798, ps20-490,bjp10-130}. In the inverse scattering transform procedure, especially in the inverse scattering problem, one important tool is the Riemann-Hilbert problem which is closely related to the sectionally meromorphic functions. In general, the contour exists self-intersection points. Thus the singularities may be eliminated by introducing the admissible contour. The solvability of the Riemann-Hilbert problem can be discussed on the the admissible contour. While the solvability of the Dbar problem can be carried out on certain domains, which is still valid even for the boundary contour with corner points (related to the self-intersection points).

In this paper, we will investigate the maps of the Dbar data relevant to the Kaup-Newell spectral problem
\begin{equation}\label{kna1}
 \begin{aligned}
&\psi_x=(-ik^2\sigma_3+k Q)\psi,  \quad Q=\left(\begin{matrix}
0&u(x)\\
v(x)&0
\end{matrix}\right),
\end{aligned}
\end{equation}
as well as to the Chen-Lee-Liu spectral problem
\begin{equation}\label{kna5}
 \begin{aligned}
\psi_x=(-ik^2\sigma_3+k Q-\frac{1}{4}iQ^2\sigma_3)\psi,
\end{aligned}
\end{equation}
and to the Gerdjikov-Ivanov spectral problem
\begin{equation}\label{kna6}
 \begin{aligned}
\psi_x=(-ik^2\sigma_3+k Q-\frac{1}{2}iQ^2\sigma_3)\psi.
\end{aligned}
\end{equation}

For any bounded domain $G\subset \mathbb{C}$, we define $C^\alpha(G)$ as the space of functions which are H\"older continuous on closure $\bar{G}$. For every element $f\in C^\alpha({G})$, its norm is defined
\begin{equation}\label{kna2}
 \|f\|_{C^\alpha({G})}=\|f\|_{L^\infty({G})}+\|f\|_{h^\alpha({G})},
\end{equation}
where
\[\|f\|_{L^\infty({G})}=\sup\limits_{k\in \bar{G}}|f(k)|<\infty
, \quad \|f\|_{h^\alpha({G})}=\sup\limits_{k_1,k_2\in \bar{G}}\frac{|f(k_1)-f(k_2)|}{|k_1-k_2|^\alpha}<\infty.\]

To apply the Dbar problem to the integrable system, we need to consider the functions defined on the whole complex plane $\mathbb{C}$. In this case, we introduce another Banach space $L^{p,\nu}(\mathbb{C})$ by $f\in L^{p,\nu}(\mathbb{C})$ admitting
\begin{equation}\label{kna3}
f\in L^p(E_1), \quad f_{(\nu)}\in L^p(E_1),
\end{equation}
where $ f_{(\nu)}(k)=|k|^{-\nu}f(k^{-1})$ and $E_1=\{k:|k|\leq1\}$. In addition,
\begin{equation}\label{kna4}
\|f\|_{L^{p,\nu}(\mathbb{C})}=\|f\|_{L^p(E_1)}+\|f_{(\nu)}\|_{L^p(E_1)}.
\end{equation}

\begin{figure}[htp]
  \centering
  \includegraphics[width=5cm]{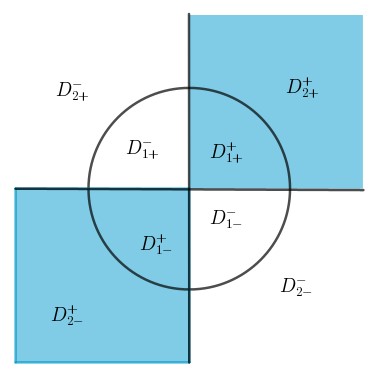}\\
  \caption{Domain $D_j^\pm$ and their split domains}\label{knfig1}
\end{figure}

For DNLS spectral problems, we need to introduce the following domains
\begin{equation}\label{knb14}
D^\pm=\{k\in\mathbb{C}:\pm{\rm Im}k^2\geq0\}.
\end{equation}
To discuss the properties of the associated integrals, we split the region of integration $D^\pm$ into two parts
\[D^\pm=D_1^\pm+D_2^\pm,\]
where
\[D_1^\pm=\{k\in D^\pm: |k|\leq1\}, \quad D_2^\pm=\{k\in D^\pm: |k|\geq1\}.\]
Moreover, the two parts $D_j^\pm, j=1,2$ are further split into $D_{j,\pm}^\pm$ (see Figure \ref{knfig1}) with
\[D_{j,\pm}^+=\{k\in D_j^+: \pm {\rm Im}k\geq0\}, \quad D_{j,\pm}^-=\{k\in D_j^-: \pm {\rm Im}k\geq0\}.\]

In addition, we need the subspaces $C^\alpha(D_{j,+}^\pm)$, $C^\alpha(D_{j,-}^\pm)$ and $C^\alpha(D^\pm)$ of $C^\alpha(\mathbb{C})$, as well as the similar subspaces of $L^{p,\nu}(\mathbb{C})$. We let $C^\alpha(\mathbb{C})=C^\alpha(D^+)\cap C^\alpha(D^-)$, where
\[ C^\alpha(D^\pm)=C^\alpha(D_1^\pm)\cap C^\alpha(D_2^\pm), \quad
C^\alpha(D_j^\pm)=C^\alpha(D_{j,+}^\pm)\cap C^\alpha(D_{j,-}^\pm).\]
The space $L^{p,\nu}(\mathbb{C})$ can be defined in a same way.

The paper is organized as follows. In section 2, we first show the properties on the integral operator $(\cdot)R\hat{T}_\mathbb{C}(k;x)$ with bounded domain. The integral operator, related to the Dbar problem without the physical variables, gives a map from $C^\alpha({G})$ to itself. Then we define the new integral operator $(\cdot)RT_\mathbb{C}(k;x)$ relevant to the DNLS spectral problems via the decomposition technique. The prior estimates are given to show the inverse operator $(\mathcal{I}-RT_\mathbb{C})^{-1}$ is exist and bounded. The associated integral equation is shown to be equivalent to a modified Dbar problem existing a unique solution. In section 3, we give the potential reconstruction of DNLS spectral problems, including the Kaup-Newell spectral problem, the Chen-Lee-Liu spectral problem and the Gerdjikov-Ivanov spectral problem, by virtue of the dressing method, as well as certain constraints and symmetry conditions. In addition, the compatibility of the modified NLS spectral problem with modified Dbar problem is discussed. Since the potential function of the Kaup-Newell spectral problem and the Chen-Lee-Liu spectral problem as well as the Gerdjikov-Ivanov spectral problem can be defined in a same way, the maps from the Dbar data $r_+\in L^{p,4}(\mathbb{C})$ to the DNLS potential $u\in L^2(\mathbb{R})$ are shown to be Lipschitz continuous.

\setcounter{equation}{0}
\section{Dbar problem relevant to the DNLS spectral problems}
\subsection{Properties on the integral operator with bounded domain}
Let $G$ be a bounded domain with continuous boundary $\Gamma$, and $\bar{G}=G\cup\Gamma$. Suppose the functions $\phi(k)\in C(\bar{G})\cap C^1(G), f(k)\in C(\bar{G})$ admitting the Dbar problem
\begin{equation}\label{knb1}
\bar\partial \phi(k)=\frac{\partial\phi(k)}{\partial\bar{k}}=f(k),
\end{equation}
then the Pompeiu formula implies that
\begin{equation}\label{knb2}
\phi(k)=\Phi(k)+f\hat{T}_G(k),
\end{equation}
where
\[\Phi(k)=\frac{1}{2\pi i}\int_{\Gamma}\frac{\phi(z)}{z-k}\mathrm{d}z,\]
\begin{equation}\label{bdab3}
 f\hat{T}_G(k)=\frac{1}{2\pi i}\iint_G\frac{f(z)}{z-k}\mathrm{d}z\wedge \mathrm{d}\bar{z}.
\end{equation}

\begin{theorem}\cite{INV1962,2603.18930v2}\label{knthb1}
Let $G$ be a bounded domain and $f\in L^p({G}), p>2$, then the function $g(k)=f\hat{T}_G(k)$ satisfying the conditions
\begin{equation}\label{knb4}
|g(k)|\leq M_1(p,G)\|f\|_{L^p({G})}, \quad k\in{G},
\end{equation}
\begin{equation}\label{knb5}
 |g(k_1)-g(k_2)|\leq M_2(p)\|f\|_{L^p(\bar{G})}|k_1-k_2|^\alpha,
\end{equation}
where $\alpha=\frac{p-2}{p}, k_1,k_2\in{G}$. Here the positive constant $M_1(p,G)$ depends on $\alpha$ and $G$, the positive constant $M_2(p)$ depends on $\mu$.
\end{theorem}

It is remarked that the Theorem \ref{knthb1}
implies that the operator $\hat{T}_G$ maps the space $L^p(\bar{G})$ to $C^\alpha(\bar{G})$, and
\begin{equation}\label{knb6}
 \|f\hat{T}_G(k)\|_{C^\alpha({G})}\leq M_3(p,G)\|f\|_{L^p({G})},
\end{equation}
where $\alpha=\frac{p-2}{p}, p>2$.

\begin{lemma}\label{knleb1}
Let $G$ be a bounded domain. If $w\in L^p({G})$ and $f\in C^\alpha({G})$, then $fw\in L^p({G})$, and
\begin{equation}\label{knb7}
\|fw\|_{L^p({G})}\leq M(p,\alpha,G)\|w\|_{L^p({G})}\|f\|_{C^\alpha({G})},
\end{equation}
where $p>1, 0<\alpha<1$.
\end{lemma}
{\bf Proof}. For every $k,k_1\in G$, let $d={\rm diam}(G)$,  we consider
\[\begin{aligned}
&\iint_G|f(k)w(k)|^p\mathrm{d}\sigma_k\\
&\leq2^{p-1}\iint_G(|f(k)-f(k_1)|^p+|f(k_1)|^p)|w(k)|^p\mathrm{d}\sigma_k\\
&\leq 2^{p-1}\left(\iint_G|f(k)-f(k_1)|^p|w(k)|^p\mathrm{d}\sigma_k+|f(k_1)|^p\iint_G|w(k)|^p\mathrm{d}\sigma_k\right)\\
&\leq 2^{p-1}\left(d^{p\alpha}\|f\|_{h^\alpha({G})}^p+\|f\|^p_{L^\infty({G})}\right)\|w\|^p_{L^p({G})}\\
&\leq \tilde{M}(p,\alpha,G)\|w\|^p_{L^p({G})}\|f\|^p_{C^\alpha({G})}.
\end{aligned}\]
This finish the proof of the lemma. \qquad $\square$

We note that $p$ and $\alpha$ in the Lemma \ref{knleb1} are mutually independent. If let $\alpha=\frac{p-2}{p}$ and $p>2$, we can obtain the following important theorem.

\begin{theorem}\label{knthb2}
Let $w\in L^{p}(G)$ and $f\in C^\alpha(\mathbb{C})$, then
\begin{equation}\label{knb8}
 \|fw\hat{T}_G(k)\|_{C^\alpha({G})}\leq M_4(p,G)\|w\|_{L^{p}({G})}\|f\|_{C^\alpha({G})},
\end{equation}
where $\alpha=\frac{p-2}{p}$, $p>2$.
\end{theorem}
{\bf Proof}. Letting $p>2$ and $\alpha=\frac{p-2}{p}$, then the condition \eqref{knb8} can be proved from \eqref{knb7} and \eqref{knb6}. \qquad $\square$

\subsection{The integral operator relevant to the DNLS spectral problems}
To investigate the DNLS spectral problems via the Dbar problem, we need to introduce the physical variable $x$ in
the spectral transformation matrix $R(k;x)$ via the following evolution equation
\begin{equation}\label{knb9}
 \partial_xR(k;x)=-ik^2[\sigma_3,R(k;x)], \quad \sigma_3={\rm diag}(1,-1),
\end{equation}
which implies that
\begin{equation}\label{knb10}
R(k;x)=\left(\begin{matrix}
0&r_+(k)\mathrm{e}^{-2ik^2x}\\
r_-(k)\mathrm{e}^{2ik^2x}&0
\end{matrix}\right).
\end{equation}

We note that the exponents $\mathrm{e}^{-2ik^2x}, x\in\mathbb{R}, k\in\mathbb{C}$ are bounded in different regions in $\mathbb{R}\times\mathbb{C}$. To control the spectral transformation matrix $R(k;x)$ in \eqref{knb10}, we split it into two nilpotent matrices $w_\pm(k;x)$, that is $R=w_-+w_+$, where
\begin{equation}\label{knb11}
w_-=\left(\begin{matrix}
0&0\\
r_-(k)\mathrm{e}^{2ik^2x}&0
\end{matrix}\right), \quad w_+=\left(\begin{matrix}
0&r_+(k)\mathrm{e}^{-2ik^2x}\\
0&0
\end{matrix}\right).
\end{equation}
Then we define a new integral operator $RT_\mathbb{C}(k;x)$ on a $2\times2$ matrix function $\psi=\psi(k;x)$ by the following decomposition
\begin{equation}\label{knb12}
\begin{aligned}
 \psi RT_\mathbb{C}(k;x):=&[\psi w_-T_{D^+}(k;x)+\psi w_+T_{D^-}(k;x)]\chi_{\{x>0\}}\\
 &+[\psi w_-T_{D^-}(k;x)+\psi w_+T_{D^+}(k;x)]\chi_{\{x<0\}},
\end{aligned}
\end{equation}
where
\begin{equation}\label{knb13}
\begin{aligned}
\psi w_-T_{D^\pm}(k;x)=\frac{1}{2\pi i}\iint_{D^\pm}\frac{\psi(z;x)w_-(z;x)}{z-k}\mathrm{d}z\wedge \mathrm{d}\bar{z},\\
\psi w_+T_{D^\pm}(k;x)=\frac{1}{2\pi i}\iint_{D^\pm}\frac{\psi(z;x)w_+(z;x)}{z-k}\mathrm{d}z\wedge \mathrm{d}\bar{z},
\end{aligned}
\end{equation}
with the domains $D^\pm$ defined in \eqref{knb14}, (see Figure \ref{knfig1}).

\begin{lemma}\label{knleb2}
Let $\psi\in L^\infty(\mathbb{R}^+,\mathbf{C}^\alpha(D^\pm))$, $r_\mp \in L_k^{p,4}(D^\pm)$, and
introduce matrix functions
\begin{equation}\label{knb15}
g^\pm(k;x):=\psi w_\mp T_{D^\pm}(k;x)\chi_{\{x>0\}},
\end{equation}
then we have the following condition
\begin{equation}\label{knb16}
\|g^\pm(\cdot;x)\|_{C^\alpha(D^\pm)}\leq M_\pm(p)\|r_\mp\|_{L^{p,4}(D^\pm)}\|\psi\|_{L^\infty(\mathbb{R}^+,C^\alpha(D^\pm))},
\end{equation}
where $\alpha=\frac{p-2}{p}, p>2$.
\end{lemma}
{\bf Proof}. we give the proof for the matrix function $g^+(k;x)$, and introduce matrix functions
\[g_{j,\pm}^+(k;x):=\psi w_-T_{D_{j,\pm}^+}(k;x)\chi_{\{x>0\}}, \quad j=1,2.\]
For the integral $g_2^\pm(k;x)$, the region of integration $D_2^\pm$ is not bounded. We consider the variable transformation $k\to k^{-1}$, then
\[D_{2,\pm}^+ \to D_{1,\mp}^-, \quad D_{2,\pm}^- \to D_{1,\mp}^+.\]
Moreover, we have \cite{2603.18930v2}
\[g_{2,\pm}^+(k;x)=g_{0,\mp}^-(0;x)-g_{0,\mp}^-(k^{-1};x), \]
where
\[g_{0,\mp}^-(\lambda;x)=\frac{1}{2\pi i}\iint_{D_{1,\mp}^-}z^2\frac{\psi_{(0)}(z;x)w_{-(4)}(z;x)}{z-\lambda}\mathrm{d}z\wedge \mathrm{d}\bar{z}.
\]
Here $f_{(\nu)}$ is defined in \eqref{kna3}.

We note that, for $x>0$, $k\in D_{1,\mp}^-$, $|k|\leq1$, and the matrix norm of $w_-(k;x)$ and $w_{-(4)}(k;x)$ satisfy the following conditions
\[\begin{aligned}
\|w_-(k;x)\|\leq|r_-(k)|, \quad k\in D^+_{1,\pm}, \\
 \|w_{-(4)}(k;x)\|\leq|r_{-(4)}(k)|, \quad k\in D^-_{1,\pm}.
\end{aligned}\]
Since $D^\pm_{1,\pm}$ are bounded domain, then, using the Theorem \ref{knthb2}, we have
\[\|g^+_{1,\pm}(\cdot;x)\|_{C^\alpha(D^+_{1,\pm})}\leq M_\pm(p,D^+_{1,\pm})\|r_-\|_{L^p(D^+_{1,\pm})}\|\psi\|_{L^\infty(\mathbb{R}^+,\mathbf{C}^\alpha(D^+_{1,\pm}))},\]
and
\[\|g^-_{0,\mp}(\cdot;x)\|_{C^\alpha(D^-_{1,\mp})}\leq \tilde{M}_\pm(p,D^-_{1,\mp})\|r_{-(4)}\|_{L^p(D^-_{1,\mp})}\|\psi_{(0)}\|_{L^\infty(\mathbb{R}^+,\mathbf{C}^\alpha(D^-_{1,\mp}))}.\]

Now we can give estimate of the matrix function $g^+(k;x)$ as
\[\begin{aligned}
g^+(k;x)=&g^+_{1,+}(k;x)+g^+_{1,-}(k;x)\\
&+g^-_{0,-}(0;x)-g^-_{0,-}(k^{-1};x)+g^-_{0,+}(0;x)-g^-_{0,+}(k^{-1};x).
\end{aligned}\]
Using the triangle inequality, we find that
\[\begin{aligned}
\|g^+(\cdot;x)\|_{C^\alpha(D^+)}&\leq M_+(p,D^+_{1,+})\|r_-\|_{L^p(D^+_{1,+})}\|\psi\|_{L^\infty(\mathbb{R}^+,\mathbf{C}^\alpha(D^+_{1,+}))}\\
&\quad+M_-(p,D^+_{1,-})\|r_-\|_{L^p(D^+_{1,-})}\|\psi\|_{L^\infty(\mathbb{R}^+,\mathbf{C}^\alpha(D^+_{1,-}))}\\
&\quad+2\tilde{M}_-(p,D^-_{1,-})\|r_{-(4)}\|_{L^p(D^-_{1,-})}\|\psi_{(0)}\|_{L^\infty(\mathbb{R}^+,\mathbf{C}^\alpha(D^-_{1,-}))}\\
&\quad+2\tilde{M}_+(p,D^-_{1,+})\|r_{-(4)}\|_{L^p(D^-_{1,+})}\|\psi_{(0)}\|_{L^\infty(\mathbb{R}^+,\mathbf{C}^\alpha(D^-_{1,+}))}\\
&\leq M_+(p,D^+_{1})\|r_-\|_{L^p(D^+_{1})}\|\psi\|_{L^\infty(\mathbb{R}^+,\mathbf{C}^\alpha(D^+_{1}))}\\
&\quad+\tilde{M}_-(p,D^-_{1})\|r_{-(4)}\|_{L^p(D^-_{1})}\|\psi_{(0)}\|_{L^\infty(\mathbb{R}^+,\mathbf{C}^\alpha(D^-_{1}))}\\
&\leq M_+(p)\|r_-\|_{L^{p,4}(D^+)}\|\psi\|_{L^\infty(\mathbb{R}^+,\mathbf{C}^\alpha(D^+))}.
\end{aligned}\]
Here we have used the notations
\[\|f\|_{L^{p,4}(D^+)}=\|f\|_{L^p(D_1^+)}+\|f_{(4)}\|_{L^p(D_1^-)},\]
\[\|f\|_{C^\alpha(D^+)}=\|f\|_{C^\alpha(D_1^+)}+\|f_{(0)}\|_{C^\alpha(D_1^-)}.\]
Thus the condition \eqref{knb16} for $g^+$ is proved. The proof of $g^-$ can be finished similarly. \qquad $\square$

\begin{lemma}\label{knleb3}
Let $\psi\in L^\infty(\mathbb{R}^+,\mathbf{C}^\alpha(D^\pm))$, $r_\pm \in L^{p,4}(D^\pm)$, and
introduce matrix functions
\begin{equation}\label{knb17}
h^\pm(k;x):=\psi w_\pm T_{D^\pm}(k;x)\chi_{\{x<0\}},
\end{equation}
then we have the following condition
\begin{equation}\label{knb18}
\|h^\pm(\cdot;x)\|_{C^\alpha(D^\pm)}\leq N_\pm(p)\|r_\pm\|_{L^{p,4}(D^\pm)}\|\psi\|_{L^\infty(\mathbb{R}^-,\mathbf{C}^\alpha(D^\pm))},
\end{equation}
where $\alpha=\frac{p-2}{p}, p>2$.
\end{lemma}
{\bf Proof}. The Lemma can be proved similarly.

\begin{theorem}\label{knthb3}
Let $\psi\in L^\infty(\mathbb{R},\mathbf{C}^\alpha(\mathbb{C}))$, $r_\pm \in L_k^{p,4}(\mathbb{C})$,
then the integral in \eqref{knb12} satisfies the following condition
\begin{equation}\label{knb19}
\|\psi RT_\mathbb{C}\|_{L^\infty(\mathbb{R},\mathbf{C}^\alpha(\mathbb{C}))}\leq M(p)(\|r_+\|_{L^{p,4}(\mathbb{C})}+\|r_-\|_{L^{p,4}(\mathbb{C})})\|\psi\|_{L^\infty(\mathbb{R},\mathbf{C}^\alpha(\mathbb{C}))},
\end{equation}
where $\alpha=\frac{p-2}{p}, p>2$.
\end{theorem}
{\bf Proof}. We note that
\[\psi RT_\mathbb{C}(k;x)=g^+(k;x)+g^-(k;x)+h^+(k;x)+h^-(k;x),\]
and collect the conditions in Lemma \ref{knleb2} and Lemma \ref{knleb3}, 
the theorem can be proved. \qquad $\square$

From the Theorem \ref{knthb3}, we know that
\begin{equation}\label{knb22}
\|RT_\mathbb{C}\|_{L^\infty(\mathbb{R},\mathbf{C}^\alpha(\mathbb{C}))\to L^\infty(\mathbb{R},\mathbf{C}^\alpha(\mathbb{C}))}\leq M(p)(\|r_+\|_{L^{p,4}(\mathbb{C})}+\|r_-\|_{L^{p,4}(\mathbb{C})}).
\end{equation}
It is remarked that the positive constant $M(p)$ depends on $p$ and the diameter 2 of the unit domain in view of the proof of Theorem \ref{knthb1} \cite{INV1962,2603.18930v2}. Generally, it cannot be required to be small.
If we let
\begin{equation}\label{knb23}
M(p)(\|r_+\|_{L^{p,4}(\mathbb{C})}+\|r_-\|_{L^{p,4}(\mathbb{C})})<1,
\end{equation}
then the inverse operator $(\mathcal{I}-RT_\mathbb{C})^{-1}$ on $L^\infty(\mathbb{R},\mathbf{C}^\alpha(\mathbb{C}))$ is exist.

\begin{corollary}\label{kncob1}
If $\psi\in L^\infty(\mathbb{R},\mathbf{C}^\alpha(\mathbb{C}))$, $r_\pm \in L^{p,4}(\mathbb{C})$, $\alpha=\frac{p-2}{p}, p>2$, then
\begin{equation}\label{knb24}
\|RT_\mathbb{C}\|_{L^\infty(\mathbb{R},\mathbf{C}^\alpha(\mathbb{C}))\to L^\infty(\mathbb{R},\mathbf{C}^\alpha(\mathbb{C}))}\leq M(p)(\|r_+\|_{L^{p,4}(\mathbb{C})}+\|r_-\|_{L^{p,4}(\mathbb{C})}).
\end{equation}
Moreover, the inverse operator $(\mathcal{I}-RT_\mathbb{C})^{-1}$ is bounded, and
\begin{equation}\label{knb25}
\begin{aligned}
\|(\mathcal{I}-RT_\mathbb{C})^{-1}\|_{L^\infty(\mathbb{R},\mathbf{C}^\alpha(\mathbb{C}))\to L^\infty(\mathbb{R},\mathbf{C}^\alpha(\mathbb{C}))}
&\leq
\frac{1}{1-M(p)(\|r_+\|_{L^{p,4}(\mathbb{C})}+\|r_-\|_{L^{p,4}(\mathbb{C})})}.
\end{aligned}
\end{equation}
\end{corollary}

In the case, there exists a unique solution $\psi\in L^\infty(\mathbb{R},\mathbf{C}_k^\alpha(\mathbb{C}))$  to the integral equation
\begin{equation}\label{knb26}
\psi(k;x)=I+\psi RT_\mathbb{C}(k;x), \quad k\in \mathbb{C}.
\end{equation}

\begin{theorem}\label{knthb4}
If $r_\pm \in L^{p,4}(\mathbb{C})$ with small norm satisfy the condition \eqref{knb23},
then there exists unique function $\psi\in {L}^\infty(\mathbb{R},\mathbf{C}^{1,\alpha}(\mathbb{C}))$ satisfying the modified Dbar equation
\begin{equation}\label{dbknb1}
\bar\partial\psi(k;x)=\left\{\begin{array}{cc}
\psi(k;x)w_-(k;x)\chi_{\{x>0\}}+\psi(k;x)w_+(k;x)\chi_{\{x<0\}}, & k\in D^+,\\
\psi(k;x)w_+(k;x)\chi_{\{x>0\}}+\psi(k;x)w_-(k;x)\chi_{\{x<0\}}, & k\in D^-,
\end{array}\right.
\end{equation}
with normalization condition $\psi(k;x)\to \psi_0(x), k\to\infty$ is equivalent to the following integral equation
\begin{equation}\label{dbknb7}
\psi(k;x)=\psi_0(x)+\psi RT_\mathbb{C}(k;x), \quad k\in \mathbb{C}.
\end{equation}
\end{theorem}
Proof. It is verified that the Dbar equation \eqref{dbknb1} with normalization condition $\psi(k;x)\to \psi_0(x), k\to\infty$ gives to the integral equation \eqref{dbknb7} by virtue of the generalized Cauchy integral formula or the Pompeiu formula
\[\phi(k)=\frac{1}{2\pi i}\int_{C^\infty}\frac{\phi(z)}{z-k}dz+\frac{1}{2\pi i}\iint_\mathbb{C}\frac{\partial_{\bar{z}}\phi(z)}{z-k}dz\wedge d\bar{z},\]
where $C^\infty$ is a circle with an infinite radius.

Conversely, for $k\in D^+$, since $\psi w_\pm T_{D^-}(k;x)$ are holomorphic,
then we have
\[\begin{aligned}
\bar\partial(\psi RT_\mathbb{C}(k;x))&=\bar\partial(\psi w_-T_{D^+}(k;x))\chi_{\{x>0\}}+\bar\partial(\psi w_+T_{D^+}(k;x))\chi_{\{x<0\}}\\
&=\psi(k;x)w_-(k;x)\chi_{\{x>0\}}+\psi(k;x)w_+(k;x)\chi_{\{x<0\}}, \quad k\in D^+.
\end{aligned}\]
Similarly, we have
\[\bar\partial(\psi RT_\mathbb{C}(k;x))=\psi(k;x)w^+(k;x)\chi_{\{x>0\}}+\psi(k;x)w^-(k;x)\chi_{\{x<0\}}, \quad k\in D^-.\]
Thus, $\psi(k;x)$ given by \eqref{dbknb7} is the solution of the Dbar equation \eqref{dbknb1}.
The theorem is proved. \qquad $\square$

\setcounter{equation}{0}
\section{Potential reconstruction of DNLS spectral problems}

Given the spectral transformation matrix $R(k;x)$ admitting equation \eqref{knb9}, we consider the piecewise Dbar problem in the Theorem \ref{knthb4} relevant to the DNLS spectral problems. If $r_\pm \in L^{p,4}(\mathbb{C})$ with small norm satisfy the condition \eqref{knb23}, then the eigenfunction $\psi(k;x)\in L^\infty(\mathbb{R},\mathbf{C}^\alpha(\mathbb{C}))$ is the unique solution of the $2\times2$ matrix integral equation
\begin{equation}\label{knc1}
\psi(k;x)=\psi_0(x)+\psi RT_\mathbb{C}(k;x), \quad k\in \mathbb{C},
\end{equation}
where the normalization $\psi_0(x)$ is a certain $2\times2$ diagonal matrix function which is related to the DNLS potentials. In addition, via the Corollary \ref{kncob1},
the unique solution can be given
\begin{equation}\label{knc2}
\psi(k;x)=\psi_0(I-RT_\mathbb{C})^{-1}(k;x).
\end{equation}

\begin{lemma}\label{dbl4}
Let $k^lr_\pm\in L^{p,4}(\mathbb{C}), (l=0,1)$  and $\psi\in L^\infty(\mathbb{R},\mathbf{C}^{\alpha}(\mathbb{C}))$, where $\psi=(\psi_{mn})_{2\times2}$. Let
\begin{equation}\label{dbd7}
\begin{aligned}
\langle f(k)\rangle_\pm=\frac{1}{2\pi i}\iint_{D^\pm}f(k)\mathrm{d}k\wedge\mathrm{d}\bar{k},\\
\end{aligned}
\end{equation}
where $D^\pm$ is defined in \eqref{knb14}, then the matrix norms admit the following conditions
\begin{equation}\label{dbd8}
|\langle k^l\psi w_\mp\rangle_\pm \chi_{\{x>0\}}|\leq M_+(p,D_1)\|\psi\|_{L^\infty(R^+,\mathbf{C}^{\alpha}(D^\pm))}\|(\cdot)^lr_\mp\|_{L^{p,4}(D^\pm)},
\end{equation}
\begin{equation}\label{dbd9}
|\langle k^l\psi k^lw_\pm\rangle_\pm\chi_{\{x<0\}}|\leq M_-(p,D_1)\|\psi\|_{L^\infty(R^-,\mathbf{C}^{\alpha}(D^\pm))}\|(\cdot)^lr_\pm\|_{L^{p,4}(D^\pm)},
\end{equation}
where $p>2$, and $\alpha=\frac{p-2}{p}$.
\end{lemma}
{\bf Proof}. We only give the proof of $\langle k^l\psi w_-\rangle_+ \chi_{\{x>0\}}$. we find that
\[\begin{aligned}
\langle k^l\psi w_-\rangle_+ \chi_{\{x>0\}}&=\frac{1}{2\pi i}\iint_{D_1^+}\psi(k;x)k^lw_-(k)\mathrm{d}k\wedge\mathrm{d}\bar{k}\chi_{\{x>0\}}\\
&\quad +\frac{1}{2\pi i}\iint_{D_2^+}\psi(k;x)k^lw_-(k)\mathrm{d}k\wedge\mathrm{d}\bar{k}\chi_{\{x>0\}}\\
&:=A_1+A_2.
\end{aligned}\]
For $A_1$, using Lemma \ref{knleb1}, we can get
\[\begin{aligned}
|A_1|&\leq\frac{1}{\pi}\iint_{D_{1+}^+}\|\psi(k;x)k^lr_-(k)\|d\sigma_k+\frac{1}{\pi}\iint_{D_{1-}^+}\|\psi(k;x)k^lr_-(k)\|d\sigma_k\\
&\leq\frac{2}{\pi}(\frac{\pi}{2})^{\frac{1}{q}}\|\psi(\cdot;x)(\cdot)^lr_-(\cdot)\|_{L^p(D_1^+)}\\
&\leq M_1(p,D_1^+)\|\psi\|_{L^\infty(R^+,\mathbf{C}^{\alpha}(D_1^+))}\|k^lr_-(k)\|_{L^p(D_1^+)},
\end{aligned}\]
where $q=\frac{p}{p-1}, p>2$. While for $A_2$, since $\psi\in L^\infty(\mathbb{R},\mathbf{C}^{\alpha}(\mathbb{C})), k^lr_-\in L^{p,4}$ then
\[\begin{aligned}
&\|\psi_{(0)}(\cdot;x)\|_{C^\alpha(E_1)}<\infty,\\
&\|(k^lr_-(k))_{(4)}\|_{L^p(E_1)}=\||k|^{-4}k^{-l}r_-(k^{-1})\|_{L^{p}(E_1)}|<\infty,
\end{aligned} \]
where $E_1=\{k:|k|\leq1\}$. Similarly, using the Lemma \ref{knleb1}, we obtain that
\[\begin{aligned}
|A_2|&\leq\frac{1}{\pi}\iint_{D_1^-}\|\psi_{(0)}(k;x)(k^lr_-(k))_{(4)}\|d\sigma_k\\
&\leq\frac{2}{\pi}(\frac{\pi}{2})^{\frac{1}{q}}\|\psi_{(0)}(\cdot;x)((\cdot)^lr_-(\cdot))_{(4)}\|_{\mathbf{L}^p(D_1^-)}\\
&\leq M_1(p,D_1^-)\|\psi_{(0)}\|_{L^\infty(R^+,\mathbf{C}^{\alpha}(D_1^-))}\|(k^lr_-(k))_{(4)}\|_{L^p(D_1^-)}.
\end{aligned}\]
Thus, using the triangle inequality, we have
\[\begin{aligned}
|\langle\psi k^lw_-\rangle_+ \chi_{\{x>0\}}|&\leq M_1(p,D_1^+)\|\psi\|_{L^\infty(R^+,\mathbf{C}^{\alpha}(D_1^+))}\|k^lr_-(k)\|_{L^p(D_1^+)}\\
&\quad+M_1(p,D_1^-)\|\psi_{(0)}\|_{L^\infty(R^+,\mathbf{C}^{\alpha}(D_1^-))}\|(k^lr_-(k))_{(4)}\|_{L^p(D_1^-)}\\
&\leq M_+(p,D_1)\|\psi\|_{L^\infty(R^+,\mathbf{C}^{\alpha}(D^+))}\|k^lr_-(k)\|_{L^{p,4}(D^+)}.
\end{aligned}\]
The Lemma is proved. \qquad $\square$

In view of conditions \eqref{dbd8} and \eqref{dbd9}, the following representations are well-defined by noting that
\begin{equation}\label{knc8}
\begin{aligned}
k\psi w_\mp T_{D^\pm}(k;x)\chi_{\{x>0\}}
=k[\psi w_\mp T_{D^\pm}(k;x)\chi_{\{x>0\}}]+\langle \psi w_\mp\rangle_\pm\chi_{\{x>0\}},\\
k\psi w_\pm T_{D^\pm}(k;x)\chi_{\{x<0\}}=k[\psi w_\pm T_{D^\pm}(k;x)\chi_{\{x>0\}}]+\langle \psi w_\pm\rangle_\pm\chi_{\{x<0\}},
\end{aligned}
\end{equation}
and
\begin{equation}\label{knc9}
\begin{aligned}
k^2\psi w_\mp T_{D^\pm}(k;x)\chi_{\{x>0\}}&=k^2[\psi w_\mp T_{D^\pm}(k;x)\chi_{\{x>0\}}]\\
&\quad +k\langle \psi w_\mp\rangle_\pm\chi_{\{x>0\}}+\langle k\psi w_\mp\rangle_\pm\chi_{\{x>0\}},\\
k^2\psi w_\pm T_{D^\pm}(k;x)\chi_{\{x<0\}}&=k^2[\psi w_\pm T_{D^\pm}(k;x)\chi_{\{x<0\}}]\\
&\quad +k\langle \psi w_\pm\rangle_\pm\chi_{\{x<0\}}+\langle k\psi w_\pm\rangle_\pm\chi_{\{x<0\}}.
\end{aligned}
\end{equation}

In terms of the definitions \eqref{knc1} and  \eqref{knb12},
equations \eqref{knc8} and \eqref{knc9} give rise to
\begin{equation}\label{knc11}
k\psi RT_\mathbb{C}=k[\psi RT_\mathbb{C}]+\langle \psi R\rangle=k\psi-k\psi_0+\langle \psi R\rangle,
\end{equation}
and
\begin{equation}\label{knc12}
\begin{aligned}
k^2\psi RT_\mathbb{C}&=k^2[\psi RT_\mathbb{C}]+k\langle \psi R\rangle+\langle k\psi R\rangle\\
&=k^2\psi-k^2\psi_0+k\langle \psi R\rangle+\langle k\psi R\rangle,
\end{aligned}
\end{equation}
where
\begin{equation}\label{knc7}
\begin{aligned}
\langle k^l\psi R\rangle=&[\langle k^l\psi w_-\rangle_+ +\langle k^l\psi w_+\rangle_-]\chi_{\{x>0\}}\\
&+[\langle k^l\psi w_+\rangle_++\langle k^l\psi w_-\rangle_-]\chi_{\{x<0\}}, \quad l=0,1,2.
\end{aligned}
\end{equation}
From \eqref{knc11}, \eqref{knc12} and using the solution formula \eqref{knc2}, we find that
\begin{equation}\label{knc16}
k\psi_0(\mathcal{I}-RT_\mathbb{C})^{-1}=k\psi+\langle \psi R\rangle\psi_0^{-1}\psi,
\end{equation}
and
\begin{equation}\label{knc17}
\begin{aligned}
k^2\psi_0(\mathcal{I}-RT_\mathbb{C})^{-1}&=k^2\psi+\langle \psi R\rangle\psi_0^{-1} k\psi_0(\mathcal{I}-RT_\mathbb{C})^{-1}\\
&=k^2\psi+k\langle \psi R\rangle\psi_0^{-1}\psi+(\langle \psi R\rangle\psi_0^{-1})^2\psi.
\end{aligned}
\end{equation}

Now differentiating of equation \eqref{knc1} with respect to $x$, and using \eqref{knc2}, \eqref{knb9}, we can get
\begin{equation}\label{knc10}
\begin{aligned}
\partial_x\psi&=(\partial_x\psi_{0}) (\mathcal{I}-RT_\mathbb{C})^{-1}\\
&\quad -ik^2\psi\sigma_3RT_\mathbb{C}(\mathcal{I}-RT_\mathbb{C})^{-1}+ik^2\psi RT_\mathbb{C}\sigma_3(\mathcal{I}-RT_\mathbb{C})^{-1}\\
&=\partial_x(\psi_{0})\psi_0^{-1}\psi
-ik^2\psi\sigma_3RT_\mathbb{C}(\mathcal{I}-RT_\mathbb{C})^{-1}+ik^2\psi RT_\mathbb{C}\sigma_3(\mathcal{I}-RT_\mathbb{C})^{-1}.
\end{aligned}
\end{equation}
For the second term on the right, since $RT_\mathbb{C}=\mathcal{I}-I(\mathcal{I}-RT_\mathbb{C})$, then \begin{equation}\label{dbd5}
RT_\mathbb{C}(\mathcal{I}-RT_\mathbb{C})^{-1}=(\mathcal{I}-RT_\mathbb{C})^{-1}-\mathcal{I}.
\end{equation}
Substituting  \eqref{dbd5}, \eqref{knc12} into \eqref{knc10}, and using the solution formula \eqref{knc2}, we can get that
\[
\begin{aligned}
\partial_x\psi&=(\partial_x\psi_{0})\psi_0^{-1}\psi-ik^2\psi\sigma_3(\mathcal{I}-RT_\mathbb{C})^{-1}+ik^2\psi\sigma_3\\
&\quad +ik^2\psi\sigma_3(\mathcal{I}-RT_\mathbb{C})^{-1}-ik^2\psi_0\sigma_3(\mathcal{I}-RT_\mathbb{C})^{-1}\\
&\quad +ik\langle \psi R\rangle\sigma_3(\mathcal{I}-RT_\mathbb{C})^{-1}+i\langle k\psi R\rangle \psi_0^{-1}\sigma_3\psi.
\end{aligned}
\]
Substitution of \eqref{knc16} and \eqref{knc17},  we have
\begin{equation}\label{knc14}
\begin{aligned}
\partial_x\psi&=(\partial_x\psi_{0})\psi_0^{-1}\psi+ik^2\psi\sigma_3+i\langle k\psi R\rangle \psi_0^{-1}\sigma_3\psi\\
&\quad -i\sigma_3[k^2\psi+k\langle \psi R\rangle\psi_0^{-1}\psi+(\langle \psi R\rangle\psi_0^{-1})^2\psi]\\
&\quad +i\langle \psi R\rangle\psi_0^{-1}\sigma_3[k\psi+\langle \psi R\rangle\psi_0^{-1}\psi]\\
&=-ik^2[\sigma_3,\psi]-ik[\sigma_3,\langle \psi R\rangle]\psi_0^{-1}\psi\\
&\quad+ \left((\partial_x\psi_{0}) -i[\sigma_3,\langle k\psi R\rangle]
 -i[\sigma_3,\langle \psi R\rangle]\psi_0^{-1}\langle \psi R\rangle\right)\psi_0^{-1}\psi.
\end{aligned}
\end{equation}

\subsection{DNLS spectral problems}
This one-order nonlinear differential equation \eqref{knc14} can be reduced to the Kaup-Newell spectral problem, the Gerdjikov-Ivanov spectral problem \cite{bjp10-130} or the Chen-Lee-Liu spectral problem \cite{ps20-490} by imposing certain constraint condition.
\begin{theorem}\label{dbth5}
If $\psi\in \mathbf{C}_x^1(\mathbb{R})\cap L^\infty(\mathbb{R},\mathbf{C}^{\alpha}(\mathbb{C}))$, and $k^lr_\pm \in L^{p,4}(\mathbb{C}), (l=0,1)$, $\alpha=\frac{p-2}{p}, p>2$,
then Kaup-Newell spectral problem can be obtained
\begin{equation}\label{dbd1}
\partial_x\psi(k;x)=-ik^2[\sigma_3,\psi(k;x)]+kQ(x)\psi(k;x),
\end{equation}
where the potential matrix $Q$ is given by the reconstruction formula about the Dbar data as following
\begin{equation}\label{dbd2}
\begin{aligned}
Q(x)\exp\left(-\frac{i}{2}\int_{-\infty}^xQ^2(y)\mathrm{d}y\sigma_3\right)=-i[\sigma_3,\langle \psi R\rangle],\\
 Q(x):=\left(\begin{matrix}
0&u(x)\\
v(x)&0
\end{matrix}\right),
\end{aligned}
\end{equation}
where $\langle \psi R\rangle$ is defined in \eqref{knc7}. In addition, the normalization function in \eqref{knc1} takes the following form
\[\psi_0(x)=\left(\frac{i}{2}\int_{-\infty}^xQ^2\mathrm{d}y\sigma_3\right).\]
\end{theorem}
{\bf Proof}. If $\psi\in \mathbf{C}_x^1(\mathbb{R})\cap L^\infty(\mathbb{R},\mathbf{C}^{\alpha}(\mathbb{C}))$, and $k^lr_\pm \in L^{p,4}(\mathbb{C}), (l=0,1)$, $\alpha=\frac{p-2}{p}, p>2$, one can obtain equation \eqref{knc14}.
By imposing the following assumptions
\begin{equation}\label{knc18}
Q=-i[\sigma_3,\langle \psi R\rangle]\psi_0^{-1},
\end{equation}
and
\begin{equation}\label{knc19}
\partial_x\psi_{0}=i[\sigma_3,\langle k\psi R\rangle]
 -Q\langle \psi R\rangle,
\end{equation}
equation \eqref{knc14} can be reduced to the Kaup-Newell spectral problem
\begin{equation}\label{knc15}
\partial_x\psi=-ik^2[\sigma_3,\psi]+kQ\psi.
\end{equation}

For convenience, we let $A^{[d]}$ and $A^{[o]}$ be the diagonal part and off-diagonal part of the matrix $A$ respectively. Equation \eqref{knc18} implies that $Q$ is an off-diagonal matrix and
\begin{equation}\label{knc20}
\langle \psi R\rangle^{[o]}=\frac{i}{2}\sigma_3Q\psi_0.
\end{equation}
Equation \eqref{knc19} can be split into
\begin{equation}\label{knc21}
\partial_x\psi_0=-Q\langle \psi R\rangle^{[o]}=\frac{i}{2}\sigma_3Q^2\psi_0,
\end{equation}
\begin{equation}\label{knc22}
 2i\sigma_3\langle k\psi R\rangle^{[o]}-Q\langle \psi R\rangle^{[d]}=0.
\end{equation}

If $Q$ takes the form in \eqref{dbd2}, then $Q^2=uv I$. From equation \eqref{knc21}, we take
\begin{equation}\label{knc23}
\psi_0=\exp\left(\frac{i}{2}\int_{-\infty}^xQ^2(y)\mathrm{d}y\sigma_3\right).
\end{equation}
The Theorem is proved. \qquad $\square$

\begin{remark}
It is remarked that the equations \eqref{knc15} and \eqref{knb9} implies that
\[\partial_x(\psi R)=-ik^2[\sigma_3,\psi R]+Qk\psi R.\]
Using the definition \eqref{knc7}, we can get that
\[\partial_x\langle \psi R\rangle=-i[\sigma_3, \langle k^2\psi R\rangle]+Q\langle k\psi R\rangle,\]
which implies that
\begin{equation}\label{knc24}
\partial_x\langle \psi R\rangle^{[d]}=Q\langle k\psi R\rangle^{[o]},
\end{equation}
by taking the diagonal part. From equations \eqref{knc24} and \eqref{knc22}, we find that
\begin{equation}\label{knc25}
 \partial_x\langle \psi R\rangle^{[d]}=\frac{i}{2}\sigma_3Q^2\langle \psi R\rangle^{[d]}.
\end{equation}
It is easy to see that $\langle \psi R\rangle^{[d]}$ and $\psi_0$ satisfy the same one-order differential equation, and we take $\langle \psi R\rangle^{[d]}=\psi_0$. Then
\[\langle \psi R\rangle=(I+\frac{i}{2}\sigma_3Q)\psi_0,\]
where $\psi_0$ is defined in \eqref{knc21}.
\end{remark}

\begin{remark}\label{knrec2}
Similarly, if put the following assumptions
\[
Q=-i[\sigma_3,\langle \psi R\rangle]\psi_0^{-1}, \eqno(\ref{knc18})
\]
and
\begin{equation}\label{knc30}
\partial_x\psi_{0}-i[\sigma_3,\langle k\psi R\rangle]
 +Q\langle \psi R\rangle=-\frac{i}{\mu}Q^2\sigma_3\psi_0,
\end{equation}
then equation \eqref{knc14} can be reduced to the following linear problem
\begin{equation}\label{knc28}
\partial_x\psi=-ik^2[\sigma_3,\psi]+kQ\psi-\frac{i}{\mu}Q^2\sigma_3\psi.
\end{equation}

Since $\psi_0$ is a diagonal matrix and $Q$ is off-diagnose, then the diagonal part of equation \eqref{knc30}
gives to
\[\partial_x\psi_0+Q\langle\psi R\rangle^{[o]}=-\frac{i}{\mu}Q^2\sigma_3\psi_0,\]
which implies that the normalization $\psi_0$ takes another form
\begin{equation}\label{knc31}
 \psi_0=\exp\left(\frac{i(\mu-2)}{2\mu}\int_{-\infty}^xQ^2(y)\mathrm{d}y\sigma_3\right).
\end{equation}

In particularly, if $\mu=2$, then the linear problem \eqref{knc28} gives the Gerdjikov-Ivanov spectral problem
\begin{equation}\label{knc32}
\partial_x\psi=-ik^2[\sigma_3,\psi]+kQ(x)\psi-\frac{i}{2}Q^2\sigma_3\psi,
\end{equation}
with the normalization $\psi_0=I$. In this case, the potential reconstruction formula will reduce to
\begin{equation}\label{knc18a}
Q=-i[\sigma_3,\langle \psi R\rangle].
\end{equation}

If let $\mu=4$, then we have the Chen-Lee-Liu spectral problem
\begin{equation}\label{knc33}
\partial_x\psi=-ik^2[\sigma_3,\psi]+kQ\psi-\frac{i}{4}Q^2\sigma_3\psi,
\end{equation}
and the associated normalization $\psi_0$ takes the same form as \eqref{knc23}.
\end{remark}

\subsection{Symmetry condition and compatibility}

In the following, we impose the constraint condition on the DNLS potential $Q$ in \eqref{knc18} or \eqref{dbd2} as following
\begin{equation}\label{knc35}
Q=H\bar{Q} H^{-1}, \quad H=\left(\begin{matrix}
0&1\\
\pm1&0
\end{matrix}\right),
\end{equation}
which implies that $v=\pm \bar{u}$ and $Q^2=\pm|u|^2 I$ are real. If $\psi_0$ defined by \eqref{knc23} is valid then $u$ admits certain boundary condition which implies that the Kaup-Newell spectral problem \eqref{dbd1} and the the Chen-Lee-Liu spectral problem \eqref{knc33} exist unique solution $\psi(k;x)$ for $k$ in $D^+$ or $D^-$. (It can be achieved by the classical spectral analysis.) Moreover, under the constraint condition \eqref{knc35}, for the Kaup-Newell spectral problem and the the Chen-Lee-Liu spectral problem, the normalization condition $\psi_0(x)$ satisfies $\psi_0(x)\overline{\psi_0(x)}=I$.
In addition, for the Kaup-Newell spectral problem \eqref{dbd1} and the the Chen-Lee-Liu spectral problem \eqref{knc33}, the associated symmetry condition for the eigenfunction $\psi(k;x)$ will be
\begin{equation}\label{knc37}
\psi(k;x)=H\overline{\psi(\bar{k};x)}H^{-1}.
\end{equation}

Now we will find the symmetry condition for the spectral transformation matrix $R(k;x)$. Noting that
\[H\sigma_3H^{-1}=-\sigma_3,  \quad (\overline{\psi_0})^{-1}H^{-1}=\psi_0H^{-1}=H^{-1}\psi_0^{-1},\]
we can get
\[\begin{aligned}
H\bar{Q}H^{-1}=-i[\sigma_3,H\overline{\langle \psi R\rangle} H^{-1}]\psi_0^{-1}.
\end{aligned}\]
which implies that
\begin{equation}\label{knc36}
\langle \psi R\rangle=H\overline{\langle \psi R\rangle} H^{-1}
\end{equation}
in view of the constraint condition \eqref{knc35} and the potential reconstruction formula \eqref{knc18}. Here $\langle \psi R\rangle$ is defined in \eqref{knc7}, that is
\[\begin{aligned}
\langle \psi R\rangle=&[\langle\psi w_-\rangle_+ +\langle \psi w_+\rangle_-]\chi_{\{x>0\}}\\
&+[\langle \psi w_+\rangle_++\langle \psi w_-\rangle_-]\chi_{\{x<0\}}.
\end{aligned}\eqno\eqref{knc7}\]

It is noted that a complex function is a map from the (original) plane $\mathbb{C}$ to the image space $\mathbb{C}$.
Taking the complex conjugate is an operation acting on the image space $\mathbb{C}$, and does not affect the measure of the original space.
Thus, from \eqref{knc7} and \eqref{knc36}, we have
\begin{equation}\label{knc38}
\langle\psi w_-\rangle_+=H\overline{\langle\psi w_+\rangle_-}H^{-1}, \quad
\langle\psi w_+\rangle_+=H\overline{\langle\psi w_-\rangle_-}H^{-1}.
\end{equation}
By the definition \eqref{dbd7}, we can get that
\[\langle\psi w_-\rangle_+=\frac{1}{2\pi i}\iint_{D^+}\psi(k)w_-(k)\mathrm{d}k\wedge\mathrm{d}\bar{k},\]
and
\[\begin{aligned}
H\overline{\langle\psi w_+\rangle_-}H^{-1}
&=\frac{-1}{2\pi i}\iint_{D^-}H\overline{\psi(k)}H^{-1}H\overline{w_-(k)}H^{-1}\mathrm{d}k\wedge\mathrm{d}\bar{k}\\
&=\frac{1}{2\pi i}\iint_{D^+}H\overline{\psi(\bar{k})}H^{-1}H\overline{w_-(\bar{k})}H^{-1}\mathrm{d}k\wedge\mathrm{d}\bar{k},
\end{aligned}\]
where integral variable substitution is used in the last equal sign.
Using the symmetry condition \eqref{knc37} and the first equation in\eqref{knc38}, we can get
\begin{equation}\label{knc39}
w_+(k)=H\overline{w_-(\bar{k})}H^{-1},
\end{equation}
which gives rise to $r_+(k)=\pm \overline{r_-(\bar{k})}$ or
\begin{equation}\label{knc40}
r_-(k)=\pm \overline{r_+(\bar{k})}.
\end{equation}
Similarly, from the symmetry condition \eqref{knc37} and the second equation in \eqref{knc38}, we also have
\[w_-(k)=H\overline{w_+(\bar{k})}H^{-1},\]
which together with \eqref{knc39} give
\begin{equation}\label{knc41}
R(k;x)=H\overline{R(\bar{k};x)}H^{-1}.
\end{equation}

For the Gerdjikov-Ivanov spectral problem \eqref{knc32} with $u\in L^2(\mathbb{R})$, we can derive the conditions \eqref{knc37} and \eqref{knc41} or \eqref{knc40} simliarly.

In other words, for a certain Lax integrable system with linear spectral problem $\partial_x\phi(x;k)=U(x;k)\phi(x;k)$, the potential matrix $U$ is an important characteristic. It is noted that the spectral transformation matrix $R(k)$ to the relevant Dbar problem will take the same role as the potential matrix $U$ to character the nonlinear integrable PDE. We have shown that the derivative NLS spectral problem is equivalent to the modified Dbar problem \eqref{dbknb1}, which exists a unique solution $\psi(k)$ if $r_\pm\in L^{p,4}(\mathbb{C})$ with the condition \eqref{knb23}.
In summary, if $r_\pm\in L^{p,4}(\mathbb{C})$ with the condition \eqref{knb23} satisfy the symmetry condition \eqref{knc40} or \eqref{knc39}, then we can find that both $H\overline{\psi(\bar{k};x)}H^{-1}$ and $\psi(k)$ solve the same modified Dbar equation \eqref{dbknb1}, which gives to symmetry condition \eqref{knc37}. That is the compatibility condition is valid.

 In the following, we will confine to consider the case for $r_+(k)\in L^{p,4}(\mathbb{C})$ and $M(p)\|r_+\|_{L^{p,4}}<1$, where $M(p)$ given in \eqref{knb22}.

\subsection{Maps for $L^2$ potential}
Since the reconstruction formula of the potential function for the Kaup-Newell spectral problem and the Chen-Lee-Liu spectral problem is given by\eqref{knc18}, the reconstruction formula for the Gerdjikov-Ivanov spectral problem is given by \eqref{knc18a}, then we can discuss them in one procedure frame.
\begin{theorem}
Let $r_+\in L^{p,4}(\mathbb{C})$, $p>2$,
then the DNLS potential $u$ satisfies the condition
\begin{equation}\label{dbd10}
|u(x)|\leq \frac{M(p)\|r_+\|_{L^{p,4}(\mathbb{C})}}{1-M(p)\|r_+\|_{L^{p,4}(\mathbb{C})}}.
\end{equation}

\end{theorem}
{\bf Proof}. For the Kaup-Newell spectral problem and the Chen-Lee-Liu spectral problem, from \eqref{dbd2}, we know that
\begin{equation}\label{knc34}
u\exp\left(\frac{i}{2}\int_{-\infty}^x|u|^2\mathrm{d}y\right)=-2i[\langle\psi_{11} r_+\rangle_+\chi_{\{x<0\}} + \langle\psi_{11} r_+\rangle_-\chi_{\{x>0\}}].
\end{equation}
Under the symmetry conditions in \eqref{knc35}, $|\psi_0|=1$. Using \eqref{dbd8}, \eqref{dbd9}, we have
\begin{equation}\label{knc27}
|u(x)|\leq M(p)\|\psi_{11}(\cdot;x)\|_{C^{\alpha}(\mathbb{C})}\|r_+\|_{L^{p,4}(\mathbb{C})}.
\end{equation}
From \eqref{knc2} and \eqref{knb25}, we can get that
\[\begin{aligned}
\|\psi(\cdot;x)\|_{C^{\alpha}(\mathbb{C})}&\leq\|(\mathcal{I}-RT_\mathbb{C})^{-1}\|_{_{C^{\alpha}(\mathbb{C})}\to_{C^{\alpha}(\mathbb{C})}}\\
&\leq \frac{1}{1-M(p)\|r_+\|_{L_k^{p,4}(\mathbb{C})}}.
\end{aligned}\]
The condition \eqref{dbd10} is proved.

For the Gerdjikov-Ivanov spectral problem, The condition \eqref{dbd10} can be  proved similarly. \qquad $\square$

In the proof of the Lemma \ref{knleb2} and Lemma \ref{knleb3}, we can confine the eigenfunction $\psi$ to the space $\psi\in L^2(\mathbb{R},\mathbf{C}^\alpha(\mathbb{C}))$. In fact, for $g^\pm$ in \eqref{knb15} and $h^\pm$ in \eqref{knb17},
using the Minkowski integral inequality, we can get that
\begin{equation}\label{knd1}
\|g^\pm\|_{L^2(\mathbb{R}^+,\mathbf{C}^\alpha(D^\pm))}\leq M_\pm(p)\|r_\mp\|_{L^{p,4}(D^\pm)}\|\psi\|_{L^2(\mathbb{R}^+,\mathbf{C}^\alpha(D^\pm))},
\end{equation}
and
\begin{equation}\label{knd1}
\|h^\pm\|_{L^2(\mathbb{R}^-,\mathbf{C}^\alpha(D^\pm))}\leq N_\pm(p)\|r_\mp\|_{L^{p,4}(D^\pm)}\|\psi\|_{L^2(\mathbb{R}^-,\mathbf{C}^\alpha(D^\pm))},
\end{equation}
where $\alpha=\frac{p-2}{p}, p>2$.

Using the symmetry condition  $r_-(k)=\pm \overline{r_+(\bar{k})}$, and from the definition \eqref{knb12} and Theorem \ref{knthb3}, we can further obtain the following conditions.
\begin{corollary}
Let $\psi\in L^2(\mathbb{R},\mathbf{C}^\alpha(\mathbb{C}))$, $r_\pm \in L^{p,4}(\mathbb{C})$,
then the integral in \eqref{knb12} satisfies the following condition
\begin{equation}\label{knd2}
\|\psi RT_\mathbb{C}\|_{L^2(\mathbb{R},\mathbf{C}^\alpha(\mathbb{C}))}\leq M(p)\|r_+\|_{L^{p,4}(\mathbb{C})}\|\psi\|_{L^2(\mathbb{R},\mathbf{C}^\alpha(\mathbb{C}))},
\end{equation}
where $\alpha=\frac{p-2}{p}, p>2$.
\end{corollary}

This corollary implies that:

\begin{corollary}
If $\psi\in L^2(\mathbb{R},\mathbf{C}^\alpha(\mathbb{C}))$, $r_\pm \in L^{p,4}(\mathbb{C})$, then
\begin{equation}\label{knd3}
\|RT_\mathbb{C}\|_{L^2(\mathbb{R},\mathbf{C}^\alpha(\mathbb{C}))\to L^2(\mathbb{R},\mathbf{C}^\alpha(\mathbb{C}))}\leq M(p)\|r_+\|_{L^{p,4}(\mathbb{C})}.
\end{equation}
Moreover, if the constraint condition \eqref{knb23} is valid, the inverse operator $(\mathcal{I}-RT_\mathbb{C})^{-1}$ exists.
\end{corollary}

\begin{corollary}
Let $r_+\in L^{p,4}(\mathbb{C})$, $p>2$,
then the DNLS potential $u$ defined satisfies the condition
\begin{equation}\label{knd6}
\|u\|_{L^2(\mathbb{R})}\leq \frac{M(p)\|r_+\|_{L^{p,4}(\mathbb{C})}}{1-M(p)\|r_+\|_{L^{p,4}(\mathbb{C})}}.
\end{equation}
In addition, the map from $r_+\in L^{p,4}(\mathbb{C})$ to $u\in L^2(\mathbb{R})$ is Lipschitz continuous.
\end{corollary}
{\bf Proof}. For the Kaup-Newell spectral problem and the Chen-Lee-Liu spectral problem, $|\psi_0|=1$, then, from \eqref{knc2}, we have
\[\|\psi\|_{L^2(\mathbb{R},\mathbf{C}^\alpha(\mathbb{C}))}\leq\frac{1}{1-M(p)\|r_+\|_{L_k^{p,4}(\mathbb{C})}},\]
and further from \eqref{knc27} we can get the condition \eqref{knd6}.

To prove the Lipschitzian continuity, we let $\tilde{r}_+(k)$ and $\tilde\psi(k;x)=(\tilde\psi_{mn}(k;x))$ are relevant to the potential $\tilde{u}(x)$ admitting the relation \eqref{knc34}, then we have
\begin{equation}\label{knd7}
\begin{aligned}
|u-\tilde{u}|&\leq |\tilde{u}| |\mathrm{e}^{i\eta}-\mathrm{e}^{i\tilde\eta}|\\
&+M(p) \big[\|\psi_{11}(\cdot;x)\|_{C^\alpha(\mathbb{C})}\|r_+-\tilde{r}_+\|_{L^{p,4}(\mathbb{C})}\\
&+M(p)\|\psi_{11}(\cdot;x)-\tilde\psi_{11}(\cdot;x)\|_{C^\alpha(\mathbb{C})}\|\tilde{r}_+\|_{L^{p,4}(\mathbb{C})}\big],
\end{aligned}
\end{equation}
in terms of conditions \eqref{dbd8} and \eqref{dbd9}, where
\[\eta=\frac{1}{2}\int_{-\infty}^\infty|u(y)|^2\mathrm{d}y,\]
and $\tilde\eta$ has a same form about $\tilde{u}$.

 Using the Lipschitzian continuity, and H\"older inequality, we have
\[\begin{aligned}
|\mathrm{e}^{i\eta}-\mathrm{e}^{i\tilde\eta}|&\leq|\eta-\tilde\eta|=\frac{1}{2}\int_{-\infty}^x\big||u|^2-|\tilde{u}|^2\big|dy\\
&\leq C_0\|u-\tilde{u}\|_{L^2(\mathbb{R})},
\end{aligned}\]
where
$C_0=\frac{1}{2}(\|u\|_{L^2(\mathbb{R})}+\|\tilde{u}\|_{L^2(\mathbb{R})})$.

Moreover, from \eqref{knc2}, we have
\[\|\psi(\cdot,x)\|_{\mathbf{C}^\alpha(\mathbb{C})}\leq\frac{1}{1-M(p)\|\tilde{r}_+\|_{L_k^{p,4}(\mathbb{C})}},\]
and further from \eqref{knc27} we can get the condition uniformly
\[|\tilde{u}|\leq \frac{M(p)\|\tilde{r}_+\|_{L^{p,4}(\mathbb{C})}}{1-M(p)\|\tilde{r}_+\|_{L^{p,4}(\mathbb{C})}}:=C_1.\]
Then we have
\[
|\tilde u|\,|e^{i\eta}-e^{i\tilde\eta}|
\le C_1 C_0\|u-\tilde u\|_{L^2(\mathbb R)}.
\]

For the term $M(p)\|\psi_{11}(\cdot;x)-\tilde\psi_{11}(\cdot;x)\|_{C^\alpha}\|\tilde r_+\|_{L^{p,4}}$ in \eqref{knd7}, we can get that
\[
M(p)\|\tilde r_+\|_{L^{p,4}}
\Big(\int_{\mathbb R}\|\psi_{11}(\cdot;x)-\tilde\psi_{11}(\cdot;x)\|_{C^\alpha}^2\,dx\Big)^{1/2}
\le M(p)\|\tilde r_+\|_{L^{p,4}}\|\psi-\tilde\psi\|_{L^2(\mathbb R,C^\alpha)}.
\]
From the relation \eqref{knc2}, we find that
\[\begin{aligned}
\psi-\tilde\psi&=I(\mathcal{I}-RT_\mathbb{C})^{-1}-I(\mathcal{I}-\tilde{R}T_\mathbb{C})^{-1}\\
&=I(\mathcal{I}-RT_\mathbb{C})^{-1}~(R-\tilde{R})T_\mathbb{C}~(\mathcal{I}-\tilde{R}T_\mathbb{C})^{-1},
\end{aligned}\]
which implies
\[\begin{aligned}
\|\psi(\cdot;x)-\tilde\psi(\cdot;x)\|_{\mathbf{C}^{\alpha}(\mathbb{C})}
&\leq \|(R-\tilde{R})T_\mathbb{C}(x,\cdot)\|_{\mathbf{C}^{\alpha}(\mathbb{C})}\\
&\quad\times(\mathcal{I}-RT_\mathbb{C})^{-1}\|_{L^\infty(\mathbb{R},\mathbf{C}^\alpha(\mathbb{C}))\to L^\infty(\mathbb{R},\mathbf{C}^\alpha(\mathbb{C}))}\\
&\quad\times(\mathcal{I}-\tilde{R}T_\mathbb{C})^{-1}\|_{L^\infty(\mathbb{R},\mathbf{C}^\alpha(\mathbb{C}))\to L^\infty(\mathbb{R},\mathbf{C}^\alpha(\mathbb{C}))}.
\end{aligned} \]
Using \eqref{knb25} and \eqref{knd3}, we have
\[\begin{aligned}
&\|\psi-\tilde\psi\|_{L^2(\mathbb{R},\mathbf{C}^{\alpha}(\mathbb{C}))}\\
&\leq\frac{M(p)\|r_+-\tilde{r}_+\|_{L^{p,4}}}{[1-M(p)(\|r_+\|_{L^{p,4}(\mathbb{C})})][1-\tilde{M}(p)(\|\tilde{r}_+\|_{L^{p,4}(\mathbb{C})})]}.
\end{aligned} \]
If let $\|r_+\|_{L^{p,4}(\mathbb{C})},\|\tilde{r}_+\|_{L^{p,4}(\mathbb{C})}\leq B$, then
\[
\|\psi-\tilde\psi\|_{L^2(\mathbb R,C^\alpha)}
\le \frac{M(p)}{(1-M(p)B)^2}\|r_+-\tilde r_+\|_{L^{p,4}}.
\]

For the term \(M(p)\|\psi_{11}(\cdot;x)\|_{C^\alpha}\|r_+-\tilde r_+\|\) in \eqref{knd7}, we also have
\[\begin{aligned}
M(p)\|r_+-\tilde r_+\|_{L^{p,4}}
\Big(\int_{\mathbb R}\|\psi_{11}(\cdot;x)\|_{C^\alpha}^2\,dx\Big)^{1/2}
&\le M(p)\|r_+-\tilde r_+\|_{L^{p,4}}\|\psi\|_{L^2(\mathbb R,C^\alpha)}\\
&\le \frac{M(p)}{1-M(p)B}\|r_+-\tilde r_+\|_{L^{p,4}}.
\end{aligned}
\]

As a result, from \eqref{knd7}, we obtain
\[\begin{aligned}
\|u-\tilde u\|_{L^2}
&\le C_1 C_0\|u-\tilde u\|_{L^2}\\
&+\frac{M(p)}{1-M(p)B}\|r_+-\tilde r_+\|_{L^{p,4}}
+\frac{M(p) B}{(1-M(p)B)^2}\|r_+-\tilde r_+\|_{L^{p,4}}.
\end{aligned}
\]
or
\[
(1-C_1 C_0)\|u-\tilde u\|_{L^2}
\le M(p,B)\|r_+-\tilde r_+\|_{L^{p,4}}.
\]
If \(C_1 C_0<1\), then we have
\[
\|u-\tilde u\|_{L^2}\le \frac{M(p,B)}{1-C_1 C_0}\|r_+-\tilde r_+\|_{L^{p,4}}.
\]


It is remarked that $u$ admits the prior estimation
\[
\|u\|_{L^2}\le \frac{M(p)B}{1-M(p)B},
\]
in view of \eqref{knd6}, then $C_0$ can be control by \(M(p)B/(1-M(p)B)\). Thus \(1-C_1 C_0\) is positive and only dependent on $p$ and \(B\) for small $B$.

For the Gerdjikov-Ivanov spectral problem, the proof of the statements can be finished similarly.
The corollary is proved. \qquad $\square$

\section*{Conclusions and Discussions}
In the paper, we extended the Dbar decomposition technique to discuss the spectral problems of the DNLS equations, including the Kaup-Newell spectral problem, the Chen-Lee-Liu spectral problem and the Gerdjikov-Ivanov spectral problem.
A modified Dbar problem was introduced and the well-posedness was discussed by virtue of a new integral operator. The DNLS potentials were reconstructed via the Dbar data, and the maps from the  Dbar data to the potentials were shown to be Lipschitz continuous.

In the discussion, the Theorem \ref{knthb1} about the integral operator $\hat{T}_G$ is of vital importance for the decomposition technique. The Pompeiu formula \eqref{knb2} is still valid for multiply-connected domains, and the boundary contours of the domains only need to be continuous. In addition, boundary contours are allowed to have corner points, which are related to the self-intersection points of contours appeared in the Riemann-Hilbert approach, and should be eliminated by considering the admissible contours. So we can say that the Dbar problem can simplify the processing procedure. Another good example is that the Dbar problem can be used to simplify long-time asymptotic analysis \cite{imrn2006-48673}.

\section*{Statements and Declarations}
\noindent{\bf Authorship contribution statements} Junyi Zhu: Conceptualization, Investigation, Methodology, Validation, and editing; Xueru Wang: Investigation, Methodology, Validation, Visualization, Writing-review; Huan Liu: Investigation, Methodology, Validation, Visualization, review.

\noindent{\bf Conflict of interest} The authors declare that they have no conflicts of interest.

\noindent{\bf Fundings}
This work was supported by National Natural Science Foundation of China (Grant No. 12571268, 12501334), Natural Science Foundation of Henan Province (Grant No. 252300421792, 262300421239), and Program for Science \& Technology Innovation Talents in Universities of Henan Province (Grant No. 26HASTIT049).


\end{document}